\documentclass[a4paper, 12pt]{article}
\usepackage{latexsym}
\usepackage{amssymb}
\newfont{\bl}{msbm10}
\newfont{\bls}{msbm7}  
\makeindex
\begin{document}
\input{amssym.def}
\input{amssym}
\newcommand{\ep}{\hspace*{\fill}$\Box$}
\newcommand{\eps}{\varepsilon}
\newcommand{\pr}{{\bf Proof. }}
\newcommand{\ms}{\medskip\\}
\newcommand{\cl}{\mbox{\rm cl}}
\newcommand{\g}{{\goth g}}
\newcommand{\gc}{${\cal G}$-complete }
\newcommand{\sa}{\stackrel{\scriptstyle s}{\approx}}
\newcommand{\prol}{\mbox{\rm pr}^{(n)}}
\newcommand{\prolo}{\mbox{\rm pr}^{(1)}}
\newcommand{\deta}{\frac{d}{d \eta}{\Big\vert}_{_{0}}}
\newcommand{\detas}{\frac{d}{d \eta}{\big\vert}_{_{0}}}
\newcommand{\R}{\mbox{{\bl R}}}
\newcommand{\N}{\mbox{{\bl N}}}
\newcommand{\C}{\mbox{{\bl C}}}
\newcommand{\Z}{\mbox{{\bl Z}}}
\newcommand{\K}{\mbox{{\bl K}}}
\newcommand{\sR}{\mbox{{\bls R}}}
\newcommand{\sN}{\mbox{{\bls N}}}
\newcommand{\gK}{{\cal K}}
\newcommand{\gR}{{\cal R}}
\newcommand{\gC}{{\cal C}}
\newcommand{\Dp}{${\cal D}'$ }                                          
\newcommand{\be}{\begin{equation}}
\newcommand{\ee}{\end{equation}}
\newcommand{\bea}{\begin{eqnarray}}
\newcommand{\eea}{\end{eqnarray}}
\newcommand{\beast}{\begin{eqnarray*}}
\newcommand{\eeast}{\end{eqnarray*}}
\newcommand{\go}{${\cal G}(\Omega)$ }
\newcommand{\grn}{${\cal G}(\R^n)$ }
\newcommand{\grp}{${\cal G}(\R^p)$ }
\newcommand{\grq}{${\cal G}(\R^q)$ }
\newcommand{\gt}{${\cal G}_\tau$ }
\newcommand{\gto}{${\cal G}_\tau(\Omega)$ }
\newcommand{\gtrn}{${\cal G}_\tau(\R^n)$ }
\newcommand{\gtrp}{${\cal G}_\tau(\R^p)$ }
\newcommand{\gtrq}{${\cal G}_\tau(\R^q)$ }
\newcommand{\gtn}{{\cal G}_\tau(\R^n) }
\newtheorem{thr}{\hspace*{-1.1mm}}[section]
\newcommand{\bt}{\begin{thr} {\bf Theorem }}
\newcommand{\et}{\end{thr}}
\newcommand{\bp}{\begin{thr} {\bf Proposition }}
\newcommand{\bc}{\begin{thr} {\bf Corollary }}
\newcommand{\blem}{\begin{thr} {\bf Lemma }}
\newcommand{\bex}{\begin{thr} {\bf Example }\rm}
\newcommand{\bexs}{\begin{thr} {\bf Examples }\rm}
\newcommand{\bd}{\begin{thr} {\bf Definition }}

\begin{center}
{\bf \Large Group Analysis of Differential Equations and\\[1mm]
 Generalized Functions}
\vskip5mm
{\large M.~Kunzinger\footnote{Supported by FWF - Research Grant
P10472-MAT of the Austrian Science Foundation.} 
and   M.~Oberguggenberger}
\end{center}
\vskip5mm


{\small  {\bf Abstract.} We present an extension of the methods 
of  classical  Lie  group analysis of differential equations to 
equations   involving  generalized  functions  (in  particular: 
distributions).  A suitable framework for such a generalization 
is  provided  by  Colombeau's theory of algebras of generalized 
functions.  We  show  that  under  some  mild conditions on the 
differential equations, symmetries of classical solutions remain 
symmetries for generalized solutions. Moreover, we introduce a 
generalization  of  the infinitesimal methods of group analysis 
that  allows  to  compute  symmetries  of  linear and nonlinear 
differential  equations  containing generalized function terms. 
Thereby, the group generators and group actions may be given by 
generalized functions themselves. 



\section{Introduction} \label{intro}  
Symmetry   properties  of  distributions  and  group  invariant 
distributional solutions (in particular: fundamental solutions)   
to   particular  types  of  linear 
differential      operators    have  been  studied  by Meth\'ee 
(\cite{Me}), Tengstrand (\cite{T}), Szmydt  and Ziemian  
(\cite{S, S2, SZ}, \cite{Z}). 
A  systematic  investigation of the transfer of classical group 
analysis   of  differential  equations  into  a  distributional 
setting is due to Berest and Ibragimov (\cite{huy,wi,ga,bi},  
\cite{amp}), again with a view to determining fundamental  
solutions of certain linear partial differential equations. A  
survey of the lastnamed studies including a comprehensive bibliography 
can be found in the third volume of \cite{crc}. 
As these approaches use methods from classical distribution  
theory,  their range is confined to linear equations and linear 
transformations of the dependent variables. 
 
Algebras of generalized functions offer the possibility of going  
beyond these limitations towards a generalization of group analysis  
to genuinely nonlinear problems involving singular terms, like  
distributions or discontinuous nonlinearities. In the present paper 
we develop a theory of group analysis of differential equations 
in algebras of generalized functions that allows a satisfactory 
treatment of such problems. This line of  
research  has been initiated in \cite{OR} and has been taken up 
in \cite{KO}. Applications to different types of algebras of  
generalized functions can be found in \cite{RR} and \cite{RW}.  
 
The  plan  of the paper is as follows: 
Section \ref{colal} provides a short introduction to the theory of 
algebras of generalized functions in the sense of J.F. Colombeau.
In section \ref{transfer} we consider 
systems   of   partial  differential  equations  together  with 
a classical symmetry group $G$ that transforms smooth solutions 
into smooth solutions. Assuming polynomial bounds on the action 
of  $G$, we can extend it to generalized functions belonging to 
Colombeau  algebras and ask whether $G$ remains a symmetry group 
for  generalized  solutions.  In section \ref{facprop} we develop methods 
based  on a factorization property of the transformed system of 
equations.   Essentially,  polynomial  bounds  on  the  factors 
suffice  to  give a positive answer. In the scalar case we show 
this  to  be  automatically  satisfied  whenever  the  equation 
contains  at least one of the derivatives of the solution as an 
isolated term. While the conditions of section \ref{facprop} concern 
some mild  assumptions  on the algebraic structure of the equations, 
section  \ref{contprop}  develops  a  topological criterion, applicable to 
systems   of  linear  equations:  the  existence  of  a  ${\cal 
C}^\infty$-continuous  homogeneous  right  inverse guarantees a 
positive  answer  as  well.  Along  the way we give examples of 
nonlinear  symmetry  transformations  of  shock  and delta wave 
solutions to linear and nonlinear systems. 
 
The  purpose  of  section \ref{ggas}  is to develop the general theory, 
allowing the equations  and  the group action (hence also its  
generators)  to  be  given  by generalized functions. Using the 
characterization  of  Colombeau  generalized functions by their 
generalized  pointvalues established in \cite{point} as well as 
results  on  Colombeau  solutions  to  ODEs,  we  show that the 
classical  procedure  for computing symmetries can be literally 
transferred to the generalized function situation. The defining 
equations  are derived as usual, but their solutions are sought 
in  generalized  functions.  This enlarges the reservoir of 
possible symmetries of classical equations and allows the study 
of  symmetries  of equations with singular terms. An example is 
provided   by   a  conservation  law  with  discontinuous  flux 
function. 
 
The  remainder  of  the  introduction  is devoted to fixing  
notations and recalling some basic definitions from group analysis 
of differential equations. 
We basically follow the 
notations and terminology of \cite{Olv}. Thus for the action of 
a  Lie  group  $G$  on some manifold $M$, assumed to be an open 
subset of some space ${\cal X} \times {\cal U}$ 
of  independent and dependent variables (with $\mbox{dim}({\cal 
X}) = p$ and $\mbox{dim}({\cal U})=q$) we write 
$g\cdot   (x,u)   =  (\Xi_g(x,u),\Phi_g(x,u))$.  Transformation 
groups  are always supposed to act regularly on $M$. If $\Xi_g$ 
does not depend on $u$, the group action is called projectable. 
Elements of 
the Lie algebra $\g$ of $G$ as well as the corresponding vector 
fields   on   $M$   will  typically be  denoted  by  ${\bf  v}$   
and  the one-parameter subgroup generated by ${\bf  v}$ by 
$\eta    \to    \exp(\eta    {\bf   v})$. $M^{(n)}$ denotes the 
$n$-jet space of $M$; the $n$-th prolongation of a 
group action $g$ or vector field ${\bf  v}$ is written as  
$\prol g$ or $\prol {\bf  v}$, respectively. 
Any system $S$ of $n$-th order differential equations 
in $p$ dependent and $q$ independent variables can be written 
in the form 
\[ 
\Delta_\nu(x,u^{(n)}) = 0, \quad 1\le \nu \le l. 
\] 
where the map 
\beast 
& \Delta: {\cal X}\times {\cal U}^{(n)} \to \R^l &\\ 
& (x,u^{(n)}) \to (\Delta_1(x,u^{(n)}),\ldots,\Delta_l(x,u^{(n)})) & 
\eeast 
will be supposed to be smooth.  
Hence $S$ is identified with the subvariety 
\[ 
S_\Delta = \{(x,u^{(n)}) : \Delta(x,u^{(n)}) = 0 \} 
\] 
of ${\cal X}\times {\cal U}^{(n)}$. For any $f:\Omega\subseteq  
{\cal X} \to {\cal U}$, $\Gamma_f$ is the graph of $f$ and 
$\Gamma_f^{(n)} := \{ (x,\mbox{\rm pr}^{(n)}f(x)) : x\in \Omega 
\}$ is the graph of the $n$-jet of $f$. 

\section{Colombeau algebras} \label{colal}
Already at a very early stage of development of the theory of 
distributions it became 
clear that it is impossible to embed the space ${\cal D}'(\Omega)$ of          
distributions over
some open subset $\Omega$ of $\R^n$ into an associative commutative algebra 
$({\cal A}(\Omega),+,\circ)$ satisfying
\begin{itemize}
\item[(i)] ${\cal D}'(\Omega)$ is linearly embedded into ${\cal A}(\Omega)$ 
and $f(x)\equiv 1$ is the unity in ${\cal A}(\Omega)$. 
\item[(ii)] There exist derivation operators $\partial_i:{\cal A}(\Omega)
\rightarrow {\cal A}(\Omega)$ ($i=1,\ldots,n$) that are linear and satisfy the 
Leibnitz rule.
\item[(iii)] $\partial_i|_{{\cal D}'(\Omega)}$ is the usual partial 
derivative ($i=1,\ldots,n$).
\item[(iv)] $\circ|_{{\cal C}(\Omega)\times {\cal C}(\Omega)}$
coincides with the pointwise product of functions. 
\end{itemize}
(Schwartz's impossibility result, \cite{Schw}). Furthermore, replacing  
${\cal C}(\Omega)$ by ${\cal C}^{(k)}(\Omega)$ does not alter this result.
On the other hand, many problems involving differentiation and nonlinearities
in the presence of singular objects require a method of coping with this
situation in a consistent manner (cf. e.g. \cite{MObook}, 
\cite{survey}, \cite{V}).
By the above, the best possible result would 
consist in constructing an algebra
${\cal A}(\Omega)$ satisfying (i)--(iii) and
\begin{itemize}
\item[(iv')] $\circ|_{{\cal C}^\infty(\Omega)\times {\cal C}^\infty(\Omega)}$
coincides with the pointwise product of functions. 
\end{itemize}
The actual construction of algebras enjoying these optimal properties 
is due to J.F. 
Colombeau (\cite{c1}, \cite{c2}, see also \cite{AB}, \cite{MObook}). 
The basic idea 
underlying his theory (in its  simplest -- the so-called  `special' -- form) 
is that of embedding the space of distributions into a factor 
algebra of ${\cal C}^\infty(\Omega)^I$
($I = (0,1]$) via regularization by convolution with a fixed  
`mollifier' $\rho \in {\cal S}(\R^n)$
with $\int \rho(x)\, dx = 1$.
In order to motivate the definition below let $\rho_\eps(x) := 
\eps^{-n}\rho(\frac{x}{\eps})$
and let $u\in {\cal E}'(\R^n)$ (the space of compactly supported distributions 
on $\R^n$). The
sequence $(u\ast \rho_\eps)_{\eps\in I}$ converges to $u$ in ${\cal D}'(\R^n)$. 
Taking this sequence as a representative of $u$ we obtain an embedding of 
${\cal D}'(\R^n)$ into the
algebra  ${\cal C}^\infty(\R)^I$. However, embedding ${\cal C}^\infty(\R^n) 
\subseteq {\cal D}'(\R^n)$ into this algebra via convolution as 
above will not yield a subalgebra since
of course $(f\ast \rho_\eps)(g\ast\rho_\eps)\not=(fg)\ast\rho_\eps$ 
in general. The idea, therefore, 
is to factor out an ideal ${\cal N}(\R^n)$ such that this 
difference vanishes in the resulting
quotient. In order to construct ${\cal N}(\R^n)$ it is obviously 
sufficient to find an
ideal containing all differences $(f\ast \rho_\eps)_{\eps\in I} - 
(f)_{\eps\in I}$.
Taylor expansion of $f\ast\rho_\eps - f$ 
shows that this term will vanish faster than any 
power of $\eps$, (uniformly on compact sets, in all derivatives) provided 
we additionally
suppose that $\int \rho(x) x^\alpha\, dx = 0$ for all $\alpha \in \N_0^n$ 
with $|\alpha|\ge 1$. The set of all such sequences
is not an ideal in ${\cal C}^\infty(\R^n)^I$, so we shall replace 
${\cal C}^\infty(\R^n)^I$
by the set of {\em moderate} sequences ${\cal E}_M(\R^n)$ whose every 
derivative is bounded 
uniformly on compact sets by some inverse power of $\eps$.

Thus we define the Colombeau algebra ${\cal G}(\Omega)$ as the quotient algebra
${\cal E}_M(\Omega)/{\cal N}(\Omega)$, where
\begin{eqnarray*}  
{\cal E}_M(\Omega)&:=&\{(u_\varepsilon)_{\eps\in I}\in 
{\cal C}^\infty(\Omega)^I: \forall K\subset\subset\Omega,  
\forall\alpha\in\N_o^n 
\mbox{ }\exists p\in \N \mbox{ with }\\&&\sup_{x\in K}|\partial^\alpha  
u_\varepsilon(x)|=O 
(\varepsilon^{-p})\mbox{ as }\varepsilon\rightarrow 0\}\\{\cal N}(\Omega)&:=& 
\{(u_\varepsilon)_{\eps\in I}\in {\cal C}^\infty(\Omega)^I:  
\forall K\subset\subset\Omega, \forall\alpha\in\N_o^n\mbox{ }\forall q\in\N\\  
&&\sup_{x\in K}|\partial^\alpha u_\varepsilon(x)|=O(\varepsilon^{q})\mbox{ as } 
\varepsilon\rightarrow 0\}. 
\end{eqnarray*} 
Equivalence  classes  of  sequences  $(u_\eps)_{\eps\in I}$  in  
${\cal G}(\Omega)$ will be denoted by $\cl[(u_\eps)_{\eps\in I}]$.  
${\cal G}(\Omega)$ is a differential algebra containing  
${\cal E}'(\Omega)$ as a linear subspace via the embedding  
$\iota:   u   \to  \cl[(u\ast\rho_\eps)_{\eps\in I}]$  depending  on  a 
mollifier $\rho\in {\cal S}(\R^n)$ as above.
$\iota$ commutes with partial derivatives 
and coincides with $u \to \cl[(u)_{\eps\in I}]$  on ${\cal D}(\Omega)$,  
thus rendering it a faithful subalgebra of ${\cal G}(\Omega)$.
The functor $\Omega \to {\cal G}(\Omega)$ is a fine sheaf of 
differential algebras
on $\R^n$ and there is a unique sheaf morphism $\hat \iota$ 
extending the above embedding
to ${\cal C}^\infty(\,.\,) \hookrightarrow {\cal D}'(\,.\,) \hookrightarrow 
{\cal G}(\,.\,)$. $\hat \iota$ commutes with partial derivatives, 
and its restriction to
$\cal C^\infty$ is a sheaf morphism of algebras.
                                                 
We shall also consider the algebra  
${\cal G}_\tau(\Omega)={\cal E}_\tau(\Omega)/{\cal N}_\tau(\Omega)$ 
of tempered generalized functions, where 
\begin{eqnarray*} 
&& {\cal O}_M(\Omega) = \{f\in {\cal C}^\infty(\Omega): 
\forall \alpha\in \N_o^n \ \exists p>0 
\ \sup\limits_{\scriptstyle x\in \Omega} (1+|x|)^{-p}  
|\partial^\alpha f(x)|<\infty\}\\ 
&& {\cal E}_\tau(\Omega) = \{ (u_\eps)_{\eps\in I}\in  
({\cal O}_M(\Omega))^I : 
\forall \alpha\in \N_o^n\ \exists p>0 \\ 
&&\qquad\qquad\quad\sup_{\scriptstyle x\in \Omega}(1+|x|) 
^{-p}|\partial^\alpha u_\eps(x)| = O(\eps^{-p})\ (\eps\to 0)\}\\ 
&& {\cal N}_\tau(\Omega) = \{ (u_\eps)_{\eps\in I}\in  
({\cal O}_M(\Omega))^I : 
\forall \alpha\in \N_o^n\ \exists p>0\ \forall \ q>0 \\ 
&&\qquad\qquad\quad\sup_{\scriptstyle x\in \Omega}(1+|x|) 
^{-p}|\partial^\alpha u_\eps(x)| = O(\eps^{q})\ (\eps\to 0)\} 
\end{eqnarray*} 
The map $\iota$ defined above is a linear embedding  
of ${\cal S}'(\R^n)$ into ${\cal G}_\tau(\R^n)$  
commuting with partial derivatives and making 
\[ 
{\cal O}_C(\R^n) = \{f\in {\cal C}^\infty(\R^n): 
\exists p>0 \ \forall \alpha\in 
\N_o^n\ \sup\limits_{\scriptstyle x\in \sR^n}  
(1+|x|)^{-p} |\partial^\alpha f(x)|<\infty\} 
\] 
a faithful subalgebra. Elements  of   ${\cal O}_M(\Omega)$  
are called {\em slowly increasing}. Componentwise insertion of  
elements of ${\cal G}$ into slowly increasing functions yields 
well defined elements of ${\cal G}$. Thus, in $\cal G$ not only polynomial
combinations of distributions (e.g. $\delta^2$) make sense but also
expressions like $\sin(\delta)$ have a well-defined meaning. 
The importance   of   ${\cal 
G}_\tau(\Omega)$  for  our  purposes  stems  from the fact that 
elements  of  this  algebra    can    even be  composed with each other 
(again by 
componentwise insertion, cf. \cite{HO}, \cite{thesis}),   
a  necessary prerequisite  for  generalizing  symmetry  methods,  
see    section  \ref{ggas}.  Especially in the theory of ODEs in the 
generalized function context it is often useful to consider the 
algebra   $\widetilde{{\cal  G}}_\tau(\Omega\times\Omega')$ 
whose    elements    satisfy    ${\cal    G}$-bounds   in   the 
$\Omega$-variables and ${\cal G}_\tau$-bounds   in   the 
$\Omega'$-variables (cf.~\cite{HO} or \cite{thesis}). 
Elements  of Colombeau algebras are usually 
denoted   by   capital  letters  with  the  understanding  that 
$(u_\eps)_{\eps\in I}$ denotes an arbitrary representative of $U \in  
{\cal    G}$. 

Nonlinear operations with distributions in ${\cal G}(\Omega)$ depend not only
on the distributions themselves but also on the regularization procedure used
in the embedding process. Thus the difference of two representatives 
$(u_\eps)_{\eps\in I}$,
$(v_\eps)_{\eps\in I}$ of generalized functions $U$ resp. $V$
may have ${\cal D}'$-limit $0$ as $\eps\to 0$ without $U$ and $V$
being equal in ${\cal G}(\Omega)$. Nevertheless $U$ and $V$ are to be 
considered
`equal in the sense of distributions' or {\em associated} with each other 
($U\approx V$).
Moreover, $U$ is called associated with some distribution $w$ 
if  $u_\eps \to w$ in ${\cal D}'$. If such a $w$ exists (which 
need not be the case, cf.
$\delta^2$), it is to be seen as the distributional
`shadow' of $U$. For example, all powers of the Heaviside function 
are associated with 
each other without being equal in the algebra itself. Also, 
$x\delta = 0$ in 
${\cal D}'(\R)$, so $x\delta \approx 0$ in ${\cal G}(\R)$,
but $x \delta \not = 0$ in ${\cal G}(\R)$. These examples 
illustrate a general principle:
assigning nonlinear properties to elements of the vector 
space ${\cal D}'(\Omega)$ 
amounts to introducing additional information which is reflected 
in a more rigid concept
of equality within ${\cal G}(\Omega)$ compared to that in ${\cal D}'(\Omega)$. 
This strict concept of equality allows for much more refined ways of 
infinitesimal modelling.  
On the ${\cal D}'$-level (the level of association) this additional 
information is lost in the
limit-process $\eps\to 0$. 

Generalized   numbers  (i.e.~the  ring  of  constants  in  case 
$\Omega$ is connected) in any of the above algebras will be  
denoted by ${\cal R}$. 
Componentwise  insertion  of  points  into  representatives  of 
generalized  functions  yields  well defined elements of ${\cal 
R}$. 

We note that there exist variants of Colombeau algebras that allow a
canonical embedding of distributions (indepenent of a fixed 
mollifier as above).
The basic idea for constructing these algebras is to replace the index set $I$ 
by the space of {\em all} possible mollifiers.
Our choice of the   special  variants  of   
Colombeau  algebras  is  aimed  at notational  simplicity.   
However,  all results presented in the sequel  carry  over   
to  the respective full variants of the algebras. Moreover, recently there have
been introduced global versions of Colombeau algebras, defined intrinsically 
on manifolds
and displaying the analogues of (i)--(iv) (with $\partial_i$ replaced 
by Lie-derivatives
with respect to smooth vector fields), see \cite{vim}. For 
applications of the theory to 
nonlinear PDEs see \cite{MObook} and the literature cited therein, 
for applications to 
mathematical physics and numerics, cf. \cite{Biag}, \cite{c3} and \cite{V}. 
\section{Transfer of Classical Symmetry Groups} \label{transfer}                
\subsection{Factorization Properties} \label{facprop}
The first question to be answered in trying to extend the applicability 
of classical group analysis to generalized solutions concerns permanence  
properties of classical symmetries: Let $G$ be the symmetry group of some  
system $S$ of PDEs and consider $S$ within the framework of  
${\cal G}(\Omega)$. Under which conditions do  
elements of $G$ also transform generalized solutions into other generalized 
solutions?  It  is the aim of this and the following section to 
answer this question. To begin with, let us fix some terminology: 
\bd  
Let $G$ be a projectable local group of transformations acting on some open  
set ${\cal M}\subseteq {\cal X}\times {\cal U}$ according to 
$g\cdot (x,u) = (\Xi_g(x),\Phi_g(x,u))$. $g$ is called slowly increasing if 
the map $u\to \Phi_g(x,u)$ is slowly increasing, uniformly for $x$ in compact  
sets. $g$ is strictly slowly increasing if $\Phi_g\in {\cal O}_M({\cal M})$.  
If $\Omega\subseteq {\cal X}$, $U\in {\cal G}(\Omega)$ and $g$ 
is (strictly) slowly increasing, the action of $g$ on $U$ is defined 
as the element 
\be \label{gafuncprgen} 
gU := \cl [((\Phi_g\circ (id\times u_\eps))\circ \Xi_g^{-1})_{\eps\in I}] 
\end{equation} 
of ${\cal G}(\Xi_g(\Omega))$. 
\et
If $U$ is a smooth function, (\ref{gafuncprgen}) reproduces the 
classical notion of group action on functions. 
Henceforth we make the tacit assumption that the differential equations 
under consideration are of a form that allows for an insertion of 
elements of Colombeau generalized functions (i.e.~the function $\Delta$ 
representing  the equations on the prolongation space is slowly 
increasing). 
Also, slowly increasing group actions are always understood to be projectable.  
Analogous to the classical setting we give the following 
\bd \label{gensg} 
Let $S$ be some system of differential equations with $p$ variables and 
$q$  unknown  functions. A solution of $S$ in  ${\cal G}$ is an 
element  $U\in  ({\cal  G}(\Omega))^q$,  with  $\Omega\subseteq 
{\cal X}$ open, which solves the system with equality in $({\cal   
G}(\Omega))^l$. A symmetry group of $S$ in ${\cal G}$ 
is a local transformation group acting on ${\cal X}\times {\cal U}$ such that  
if $U$ is a solution of the system in ${\cal G}$, $g\in G$ and $g\cdot U$ is  
defined, then also $g\cdot U$ is a solution of $S$ in ${\cal G}$. 
\et 
Let us take a look at the transition problem from classical to generalized 
symmetry groups on the level of representatives. Thus, let $G$ be a  
slowly increasing symmetry group of some differential 
equation 
\be \label{diffe} 
\Delta(x,u^{(n)}) = 0. 
\end{equation} 
This means that if $f$ is a classical solution, i.e.~if  
$\Delta(x,\prol f (x)) = 0$ for all $x$ then also 
$\Delta(x,\prol (g\cdot f) (x)) = 0$. Now let $U\in {\cal G}(\Omega)$ be 
a generalized solution to (\ref{diffe}). Then for any representative 
$(u_\eps)_{\eps\in I}$ of $U$ there exists some $(n_\eps)_{\eps\in I} \in  
{\cal N}(\Omega)$ such that for all $x$ and all $\eps$ we have 
\be \label{differep} 
\Delta(x,\prol u_\eps(x)) = n_\eps(x). 
\end{equation} 
In particular, the differential equation (\ref{diffe}) need not be 
satisfied for even one single value of $\eps$. This basic observation 
displays quite fundamental obstacles to a direct 
utilization of the classical symmetry group properties of $G$ in order 
to obtain statements on the status of $G$ in the Colombeau-setting.  
Therefore we have to derive properties of symmetry groups that are  
better suited to allow a transfer to differential algebras. 
The   starting   point  for  our  considerations  is  a  slight 
modification  of  a well known factorization property of smooth 
maps (cf.~\cite{Olv} , Proposition 2.10):  
\bp \label{olverprop} 
Let $F$ be a smooth mapping from some manifold $M$ to $\R^k$ 
($k\le n=\mbox{dim}(M)$), let $f: (-\eta_o,\eta_o) \times M 
\to \R$ be smooth and suppose that $f(\eta,\,.\,)$ vanishes on 
the zero set ${\cal S}_F$ of $F$, identically in $\eta$. 
If $F$ is of maximal rank ($=k$) on ${\cal S}_F$ 
then there exist smooth functions $Q_1,\ldots,Q_k: 
(-\eta_o,\eta_o)\times M \rightarrow \R$ such that  
\[ f(\eta,m)=Q_1(\eta,m)F_1(m)+\ldots+Q_k(\eta,m)F_k(m) \] 
for all $(\eta,m)\in (-\eta_o,\eta_o)\times M$. 
\ep
\et 

We  are  mainly  interested  in  the  following  application of 
Proposition \ref{olverprop}: 
\bt 
\label{mainfactor} 
Let  
\be \label{system} 
\Delta_\nu(x,u^{(n)}) = 0, \quad 1\le \nu \le l 
\end{equation} 
be a nondegenerate system of PDEs. Let $G=\{g_{\eta}:\eta \in  
(-\eta_o,\eta_o)\}$ be a one parameter symmetry group of  
(\ref{system}) and set 
$g_\eta\cdot (x,u) = (\Xi_\eta(x,u),\Phi_\eta(x,u))$. 
Then there exist  
${\cal C}^\infty$-functions $Q_{\mu \nu}: 
(-\eta_o,\eta_o)\times {\cal V} \rightarrow \R$ ($1\le \mu,\nu \le l$, 
${\cal V}$ an open subset of ${\cal M}^{(n)}$)  
such that if $u:\Omega\subseteq \R^p\rightarrow \R^q$ is smooth 
and $g_\eta u$ exists we have 
\bea 
& \Delta_{\nu}(\Xi_{\eta}(x,u(x)),\prol(g_\eta u)(\Xi_{\eta}(x,u(x))))= &  
\nonumber\\ 
& = \sum\limits_{\mu=1}^l Q_{\mu \nu}(\eta,x,\prol u(x)) 
\Delta_{\mu}(x,\prol u(x)) \label{glg4} & 
\eea 
on the domain of $g_\eta u$ for $1 \le \nu \le l$. 
\et 

\pr Denote by $z$ the coordinates on ${\cal M}^{(n)}$. That $g_\eta$ is an  
element of the symmetry group of the system is equivalent with 
\[ 
\Delta(z)=0\,\Rightarrow \, 
\Delta_\nu(\prol g_\eta(z))=0 \quad (1\le \nu \le l) 
\] 
for all $\eta$ and $z$ such that this is defined. 
$\Delta$ is of maximal rank because (\ref{system}) is nondegenerate. 
Hence, by Proposition \ref{olverprop} there exist 
${\cal C}^\infty$-functions $Q_{\mu \nu}:(-\eta_o,\eta_o)\times  
{\cal V}  \rightarrow \R$ ($1\le \mu \le l$, ${\cal V}$ an open subset of 
${\cal M}^{(n)}$) such that  
\be \label{glg5} 
\Delta_{\nu}(\prol g_\eta(z))= 
\sum_{\mu=1}^l Q_{\mu \nu}(\eta,z)\Delta_{\mu}(z). 
\end{equation} 
Now for a smooth function $u:\Omega\subseteq\R^p\rightarrow  
\R^q$ as in our assumption and $x\in \Omega$ we set 
\be \label{glg5a} 
z_u(x):=(x,\prol u(x))\in {\cal M}^{(n)}. 
\end{equation} 
Then by definition $\prol g_\eta(z_u(x))= 
(\Xi_{\eta}(x,u(x)),\prol (g_\eta u)(\Xi_{\eta}(x,u(x))))$, so the result 
follows.\ep\ms 
For a single PDE $\Delta(x,\prol u)=0$, equation (\ref{glg4}) takes the  
simpler form 
\be 
\label{glg7} 
\quad \Delta(\Xi_{\eta}(x,u(x)),\prol (g_\eta u)(\Xi_{\eta}(x,u(x))))= 
Q(\eta,x,\prol u(x))\Delta(x,\prol u(x)). 
\end{equation} 
Theorem \ref{mainfactor} will be one of our main tools in transferring  
classical symmetry groups of (systems of) PDEs into the setting of algebras of  
generalized functions. 
\bp \label{emfactor} 
Let $\eta\rightarrow g_\eta$ be a  
slowly increasing one parameter symmetry group of (\ref{system}).  
If $P_{\mu \nu}:=(Q_{\mu \nu}(\eta,\Xi_{-\eta}(\,.\,), 
\prol u_\eps(\Xi_{-\eta}(\,.\,))))_{\eps\in I}$ belongs to  
${\cal E}_M(\Omega)$ 
for $1\le \mu,\nu\le l$ and every $(u_\eps)_{\eps\in I}\in {\cal E}_M(\Omega)$, 
then  $\eta\rightarrow g_\eta$ is a symmetry group of (\ref{system})  
in ${\cal G}$ as well. This condition is satisfied if 
\[ 
(x,u^{(n)}) \to Q_{\mu\nu}(\eta,x,u^{(n)})  
\] 
is slowly increasing in the $u^{(n)}$-variables, uniformly in $x$ on 
compact sets for $1\le \mu,\nu\le l$ and every $\eta$. 
\et 

\pr It suffices to observe that (\ref{glg4}) gives 
\beast 
&& \Delta_{\nu}(x,\prol (g_\eta u)(x))=\\ 
&& = \sum_{\mu=1}^l Q_{\mu \nu}(\eta,\Xi_{-\eta}(x), 
\prol u(\Xi_{-\eta}(x))) 
\Delta_{\mu}(\Xi_{-\eta}(x),\prol u(\Xi_{-\eta}(x))). 
\eeast 
For any solution $U\in {\cal G}(\Omega)$ with representative  
$u=(u_\eps)_{\eps\in I}$, 
this expression is in ${\cal N}(\Omega)$ since $P_{\mu\nu}\in 
{\cal E}_M(\Omega)$ for each $\mu,\nu$, and every  
$\Delta_{\mu}(\Xi_{-\eta}(\,.\,),\prol u(\Xi_{-\eta}(\,.\,)))$ is 
in ${\cal N}(\Omega)$ because $U$ is a solution   
and $\Xi_{-\eta}$ is a diffeomorphism.\ep 
\bex \label{noncon} 
The system 
\begin{eqnarray} 
U_t+UU_x &=& 0 \label{s}\nonumber\\ 
V_t+UV_x &=& 0\\     
U\mid_{\{t=0\}}=U_o &,& V\mid_{\{t=0\}}=V_o\nonumber 
\end{eqnarray} 
may serve as a simplified model for a one-dimensional, elastic
material of high density in a nearly plastic state. It was
analyzed in \cite{case}, where solutions
$U,V\in {\cal G}_{s,g}(\R\times [0,\infty))$, $U_o,V_o\in 
{\cal    G}_{s,g}(\R)$ were constructed and studied (${\cal G}_{s,g}$ is a 
variant  of  the Colombeau algebra with global instead of local 
bounds). In the following we present 
some  applications  of  the above results to this system (for a 
more detailed study, see \cite{KO}). 
For $U_o'\ge 0$ (\ref{s}) 
has a unique solution $(U,V)$ in ${\cal G}_{s,g}(\R\times [0,\infty))$ with 
$\partial_x U\ge 0$. We consider solutions in 
${\cal G}_{s,g}(\R\times [0,\infty))$ with initial data  
$U_o(x)=u_{\scriptscriptstyle L}+ 
(u_{\scriptscriptstyle R}-u_{\scriptscriptstyle L})H(x)$ and 
$V_o(x)=v_{\scriptscriptstyle L}+ 
(v_{\scriptscriptstyle R}-v_{\scriptscriptstyle L})H(x)$,  
where $H$ is a generalized Heaviside function with $H'\ge 0$,  
i.e.    $H$    is    a  member  of  ${\cal G}_{s,g}(\R)$ with a 
representative  $(h_\eps)_{\eps\in I}$  coinciding  with the classical 
Heaviside function $Y$ off the interval $[-\eps,\eps]$.  
For $u_{\scriptscriptstyle L}<u_{\scriptscriptstyle R}$  
the solution $(U,V)$ is associated with the rarefaction wave 
\be 
\label{ufan} 
u(x,t)=\left\{\begin{array}{l@{,\qquad}l}u_{\scriptscriptstyle L} &  
x\le u_{\scriptscriptstyle L}t\\  
\frac{x}{t} & u_{\scriptscriptstyle L}t\le x\le  
u_{\scriptscriptstyle R}t\\ u_{\scriptscriptstyle R} &  
u_{\scriptscriptstyle R}t\le x\end{array}\right. 
\end{equation} 
\be 
\label{vfan} 
v(x,t)=\left\{\begin{array}{c@{,\qquad}l}v_{\scriptscriptstyle L} &  
x\le u_{\scriptscriptstyle L}t\\  
\left(\frac{v_{\scriptscriptstyle R}- 
v_{\mbox{\tiny L}}}{u_{\scriptscriptstyle R} 
-u_{\mbox{\tiny L}}}\right)\frac{x}{t}+\left( 
\frac{v_{\mbox{\tiny L}}u_{\scriptscriptstyle R}- 
v_{\scriptscriptstyle R}u_{\mbox{\tiny L}}}{u_{\scriptscriptstyle R}- 
u_{\mbox{\tiny L}}}\right)  
& u_{\scriptscriptstyle L}t\le x\le u_{\scriptscriptstyle R}t\\  
v_{\scriptscriptstyle R} & u_{\scriptscriptstyle R}t\le x\end{array}\right. 
\end{equation} 
However, choosing different generalized Heaviside functions for modelling 
the initial data $U_o$, respectively $V_o$ we may obtain a 
superposition of the rarefaction wave (\ref{ufan}) in $u$ with a shock 
wave 
\be \label{vsup} 
v(x,t) = v_{\scriptscriptstyle L}+(v_{\scriptscriptstyle R}  
    - v_{\scriptscriptstyle L})Y(x-ct)  
\end{equation} 
with arbitrary shock speed $c$, $u_{\scriptscriptstyle L}\le c\le 
u_{\scriptscriptstyle  R}$.  We  are  going  to construct a one 
parameter   symmetry   $\eta\to   g_\eta$  of  (\ref{s})  which 
transforms any of the solutions (\ref{vfan}), (\ref{vsup}) into 
a shock wave solution as $\eta\to\pm \infty$. For this we employ the 
two-dimensional Lorentz-transformation $(\eta,(x,t))\rightarrow  
(x\cosh(\eta)-t\sinh(\eta),-x\sinh(\eta)+t\cosh(\eta))$ with infinitesimal 
generator $X_o= -t\partial_x-x\partial_t$. Then  
$X:=X_o+(u^2-1)\partial_u$ generates a projectable one-parameter symmetry  
group  of  (\ref{s}).  Assuming that $-1 <u_{\mbox{\tiny L}} < 
u_{\mbox{\tiny  R}} <1$, we can extend the solution $(U,V)$ to 
the   region   $\Omega   =  \R^2\setminus  \{(x,t)  :  t\le  0, 
u_{\mbox{\tiny  R}}t\le  x  \le u_{\mbox{\tiny L}}t\}$ by the 
method  of characteristics applied to representatives. Then the 
Lorentz-transformed solutions  
\begin{eqnarray} 
\tilde{u}_\varepsilon(x,t) & = & -\tanh(\eta - \mbox{Artanh}(u_\varepsilon 
(x\cosh(\eta)+t\sinh(\eta),\nonumber\\ 
& & x\sinh(\eta)+t\cosh(\eta))))\label{utrans}\\ 
\tilde{v}_\varepsilon(x,t) & = & v_\varepsilon 
(x\cosh(\eta)+t\sinh(\eta),x\sinh(\eta)+t\cosh(\eta))\label{vtrans} 
\end{eqnarray} 
(with $\mbox{Artanh}$ the inverse of $\mbox{tanh}$) are well defined at least 
on $\R\times (0,\infty)$.  
The factorization property (\ref{glg4}) in this case reads   
\begin{eqnarray} 
\lefteqn{ 
(\partial_t \tilde{u}_\varepsilon + \tilde{u}_\varepsilon \partial_x 
\tilde{u}_\varepsilon)(x,t) =}  
\label{fact}\\ 
&\hspace{-0.8cm}\left((\partial_t u_\varepsilon +  
u_\varepsilon  
\partial_x u_\varepsilon)/(\cosh^3(\mbox{Artanh}(u_\varepsilon-\eta)) 
\cosh(\mbox{Artanh}(u_\varepsilon)))\right) 
\left(\Xi_\eta^{-1}(x,t)\right) \nonumber 
\end{eqnarray} 
and   similarly   for   the   second   line   in  (\ref{s}), 
demonstrating  that  $(\widetilde{U},\widetilde{V})$ is again a 
solution. For each $\eta$, $\widetilde{U}$ is associated with a 
piecewise smooth function which converges to $\mp 1$ as  
$\eta\to \pm \infty$. Observing that the coordinate transformations    
in (\ref{utrans}),   (\ref{vtrans}) approach  boosts   in   the 
directions $(\mp 1,1)$ as $\eta\to \pm \infty$, we see that the 
functions associated with $\widetilde{V}$ converge to the shock 
wave $v_{\scriptscriptstyle L}+(v_{\scriptscriptstyle R}  
- v_{\scriptscriptstyle L})Y(x\pm t)$ as $\eta\to \pm \infty$, 
for   whatever   solution   $V$   given   in   (\ref{vfan})  or 
(\ref{vsup}). 
\et 
Although   Proposition \ref{emfactor}  provides  a  manageable  
algorithm to determine if classical symmetry groups carry over to 
generalized solutions it  would certainly be preferable to have  
criteria at hand that allow to judge directly from the given  
PDE if the factors $P_{\mu\nu}$ behave nicely (given slowly  
increasing group actions). 
The first step in this direction is gaining control over the behaviour 
of the map $z\rightarrow \prol g_{\eta}(z)$, defined on ${\cal M}^{(n)}$. 
\bp \label{prolprop} 
If $\eta\rightarrow g_\eta$ is a (strictly) slowly increasing group action 
on ${\cal M}$ then $z\rightarrow \prol g_{\eta}(z)$ is (strictly) 
slowly increasing as well. 
\et 

\pr 
Let $N:=\mbox{dim}({\cal M}^{(n)})$. For 
$z=(z_1,\ldots,z_p,z_{p+1},\ldots,z_q,\ldots,z_N)\in {\cal M}^{(n)}$  
we choose some smooth function $h:{\cal X}\rightarrow {\cal U}$ satisfying 
$z=z_h(z_1,\ldots,z_p)$,   with $z_h(x)$ as in (\ref{glg5a}).  
Then we set $x:=(z_1,\ldots,z_p)$,  
$u=(z_{p+1},\ldots,z_q)$, $\widetilde{x}=\Xi_\eta(x)$  
and $\widetilde{u}=\Phi_\eta(x,u)$. By the definition of prolongued  
group actions we have to find estimates for every 
\be 
\label{glg9} 
A_s:=\left((\Phi_\eta\circ (id\times h))\circ \Xi_{-\eta}  
\right)^{(s)}(\widetilde{x}) 
\end{equation} 
(where $(s)$ denotes the derivative of order $s$) in terms of $z$. 
The above formula contains the components of $\prol g(z)$ of order $s$ 
($s\le n$).  
Note that the   
particular choice of $h$ has no influence on (\ref{glg9}), i.e.~$A_s$ 
depends exclusively on $z$. 
To compute $A_s$ explicitly we use the formula for higher order derivatives 
of composite functions (see \cite{Fr}). 
Denoting by $\Upsilon_m$ the group of permutations of $\{1,\ldots,m\}$ we have: 
\be 
\label{glg10} 
A_s(r_1,\ldots,r_s)=\sum_{i=1}^s\sum_{k\in \sN^{i} \atop |k|=s} 
\sum_{\sigma\in \Upsilon_s}\frac{1}{i!k!}(\Phi_\eta\circ (id\times h))^{(i)} 
((\widetilde{x}))(t_1,\ldots,t_i), 
\end{equation} 
where 
\[ 
t_1=\Xi_{-\eta}^{(k_1)} (\widetilde{x}) 
(r_{\sigma(1)} ,\ldots,r_{\sigma(k_1)}),\hspace{1mm}. 
\hspace{1mm}.\hspace{1mm}. 
\hspace{1mm},t_i=\Xi_{-\eta}^{(k_i)} (\widetilde{x}) 
  (r_{\sigma(s-k_i+1)} ,\ldots,r_{\sigma(s)}). 
\] 
and 
\be 
\label{glg11} 
\hspace*{1em} (\left((\Phi_\eta\circ 
(id\times h))\right)^{(i)}(x)(t_1,\ldots,t_i)= 
\sum_{j=1}^i\sum_{l\in \sN^{j} \atop |l|=i} 
\sum_{\tau\in \Upsilon_i}\frac{1}{j!l!}\Phi_\eta^{(j)}(x,u)(s_1,\ldots,s_j), 
\end{equation} 
where 
\[ 
s_1=(id\times h)^{(l_1)}(x)(t_{\tau(1)} ,\ldots,t_{\tau(l_1)}), 
\hspace{1mm}.\hspace{1mm}.\hspace{1mm}. 
\hspace{1mm},s_j=(id\times h)^{(l_j)}(x)(t_{\tau(i-l_j+1)} , 
\ldots,t_{\tau(i)}). 
\] 
Each $s_m$ consists of sums of products of certain $t_{\tau(k)}$ with  
certain $z_l$ and an analogous assertion holds for the  
$\Phi_\eta^{(j)}(x,u)(s_1,\ldots,s_j)$.  
Hence from (\ref{glg10}) and (\ref{glg11}) the result follows.\ep\ms 
Returning to our original task of finding a priori estimates for the  
factors $P_{\mu\nu}$, even with the aid of Proposition \ref{prolprop}  
we still need some information about the explicit form of the $Q_{\mu\nu}$ 
to go on. In general this seems quite difficult to achieve.  
However, there is a large and important 
class of PDEs that allow a priori statements on the concrete form of the 
factorization.  Namely,  we  are going to show that each scalar 
PDE  in which at least $u$ or one of its derivatives appears as 
a single term with constant coefficient belongs to this class. 
 
Consider a scalar PDE $\Delta(x,u^{(n)})=0$  
together with a symmetry group $\eta\rightarrow g_\eta$.  
Then we have  
\[ 
\Delta(z)=0\quad \Rightarrow \quad 
\Delta(\prol g_\eta(z))=0 
\] 
Set $F(z):=\Delta(z)$, $f(z):=\Delta(\prol g_\eta(z))$ 
and $N=\mbox{dim}({\cal M}^{(n)})$. 
Suppose that in a neighborhood of some $\bar{z}$ with  
$F(\bar{z})=0$ we have $\frac{\partial F}{\partial z_k}>0$ for some  
$1\le k\le N$.  
Then by the implicit function theorem, locally there exists a smooth function 
$\psi:\R^{N-1}\rightarrow\R$ such that in a suitable neighborhood of  
$\bar{z}$ we have 
\[F(z)=0\quad \Leftrightarrow  \quad z_k=\psi(z'),\] 
where $z'=(z_1,\ldots,\hat{z_k},\ldots,z_N)$ (meaning that the component $z_k$ 
is missing from $z'$). 
It follows that 
\[ 
F(z) = (z_k-\psi(z'))\int_0^1 \frac{\partial F}{\partial z_k} 
(z_1,\ldots,z_{k-1},\tau z_k+(1-\tau)\psi(z'),\ldots,z_N)\,d\tau, 
\] 
and on the other hand 
\[ 
f(z) = (z_k-\psi(z'))\int_0^1 \frac{\partial f}{\partial z_k} 
(z_1,\ldots,z_{k-1},\tau z_k+(1-\tau)\psi(z'),\ldots,z_N)\,d\tau.  
\] 
Thus in the said neighborhood we have 
\be 
\label{glg12} 
f(z)=F(z)\frac{\int_0^1 \frac{\partial f}{\partial z_k}(z_1,\ldots,z_{k-1}, 
\tau z_k+(1-\tau)\psi(z'),\ldots,z_N)d\tau} 
{\int_0^1 \frac{\partial F}{\partial z_k}(z_1,\ldots,z_{k-1}, 
\tau z_k+(1-\tau)\psi(z'),\ldots,z_N)d\tau} 
\end{equation} 
provided the denominator of this expression is $\not=0$. 
In particular, if for some constant $c\not=0$ we have 
$\frac{\partial F}{\partial z_k}\equiv c$ in a neighborhood of $\bar{z}$ 
then (\ref{glg12}) simplifies to 
\be 
\label{glg13} 
f(z)=\frac{1}{c}F(z)\int_0^1 \frac{\partial f}{\partial z_k} 
(z_1,\ldots,z_{k-1}, 
\tau z_k+(1-\tau)\psi(z'),\ldots,z_N)d\tau 
\end{equation} 
After these preparations we can state 
\bt 
\label{ogth} 
Let  
$\eta\rightarrow g_\eta$ be a slowly increasing symmetry group of  
the equation $\Delta(x,u^{(n)})=0$.  
Set $N=\mbox{dim}({\cal M}^{(n)})$ and suppose that   
$\frac{\partial \Delta}{\partial z_k}\equiv c\not=0$ for  
some $p+1\le k\le N$.  
Then $\eta\rightarrow g_\eta$ is a symmetry group of  
$\Delta(x,u^{(n)})=0$ in ${\cal G}$. 
\et

\pr Without loss of generality we may assume $c=1$. 
Using the above notations we have $F(z)=z_k-\psi(z')$, so (\ref{glg13}) 
implies 
\[ 
f(z)=F(z)\int\limits_0^1 \frac{\partial f}{\partial z_k}(z_1,\ldots,z_{k-1}, 
\tau z_k+(1-\tau)(z_k-F(z)),\ldots,z_N)d\tau =:F(z)Q(\eta,z). 
\] From Proposition \ref{prolprop} we know that $z\rightarrow f(z)$ is slowly 
increasing in the $u^{(n)}$-variables (i.e. in those $z_i$ with $i>p$), 
uniformly in $x=(z_1,\ldots,z_p)$ on compact sets. 
Since $F$ is slowly increasing we infer that $Q(\eta,z_u(x))\in  
{\cal E}_M(\Omega)$ for any $u\in {\cal E}_M(\Omega)$ 
(with $z_u$ as in (\ref{glg5a})). Finally, 
\[ 
\Delta(x,\prol (g_\eta u)(x))= 
\Delta(\Xi_{-\eta}(x),\prol u(\Xi_{-\eta}(x))) 
Q(\eta,\Xi_{-\eta}(x),\prol u(\Xi_{-\eta}(x))). 
\] 
Since $\Xi_{-\eta}$ is a diffeomorphism, it follows that if  
$U=\cl[u]$ solves the equation, so does $g_\eta U$.\ep\ms 
As the proof shows, we can drop the assumption $p+1\le k$ 
if we require the group action to be strictly slowly increasing. 
It is clear that many PDEs satisfy the requirements of Theorem \ref{ogth}. 
For example, in the Hopf equation $\Delta(x,t,u,u_x,u_t)=u_t+uu_x$  
or  $\Delta(z_1,\ldots,z_5)=z_5+z_3  z_4$  one  can take $k=5$. 
Note however that not every 
symmetry   group  of  this  equation  is  automatically  slowly 
increasing. Theorem \ref{ogth} constitutes a useful tool for transferring 
classical symmetry groups to Colombeau algebras. 
\bex \label{grex0} 
We consider the initial value problem for the nonlinear transport
equation 
\be \label{semi} 
\begin{array}{l} 
U_t + \lambda\cdot \nabla_x U = f(U) \\   
U\mid_{\{t=0\}}=U_o 
\end{array} 
\end{equation} 
with $t \in \R, x, \lambda \in \R^n$. It has unique solutions in 
${\cal G}(\R^{n+1})$, given $U_o \in {\cal G}(\R^n)$, 
provided   $f\in   {\cal   O}_M$  is  globally  Lipschitz  (see 
\cite{MObook}).  If  in addition $f$ is bounded and the initial 
data are distributions with discrete support, say  
$U_0(x)  =  \sum_{i,j}  a_{ij}  \delta^{(i)}(x-\xi_j)$ with
$\xi_j \in \R^n, i \in \N_0^n$, then the 
generalized  solution  is  associated with the delta wave $v+w$ 
where   
\be \label{dw} 
v(x,t)  =  \sum_{i,j}  a_{ij}  \delta^{(i)}(x-\lambda t 
-\xi_j) 
\end{equation}    
and $w$ is the smooth solution to $w_t+\lambda\cdot\nabla_x w=f(w)$, $w(0)=0$. 

The vector field $X = cf(u)\partial_u$ generates an infinitesimal 
symmetry of (\ref{semi}) for arbitrary $c \in \R$. With 
$F(u) := \int\,du/f(u)$, the corresponding Lie point transformation is 
\be \label{trans} 
(x,t,u) \to  
(\widetilde{x},\widetilde{t},\widetilde{u}) =  
(x,t,F^{-1}\left( c\eta + F(u)\right)). 
\end{equation} 
This provides a well-defined nonlinear transformation of  
the generalized solution $U \in {\cal G}(\R^{n+1})$, provided that  
the right hand side in (\ref{trans}) is slowly increasing.  
 
In the example 
\be \label{grex01} 
     U_t + \lambda\cdot\nabla_x U = \tanh (U) 
\end{equation} 
the  generalized  solution  is associated with $v(x,t)$ and $w$ 
vanishes  identically. Applying (\ref{trans}) we obtain (due to 
Theorem \ref{ogth}) the new generalized solution 
\be \label{grex02} 
  \widetilde{U}(x,t) =  
  \mbox{Arsinh}\left( e^{c\eta}\sinh(U(x,t))\right) 
\end{equation} 
(with $\mbox{Arsinh}$ the inverse of $\sinh$).
We  are  going to show that $\widetilde{U}$ is still associated 
with the delta wave $v$ in (\ref{dw}). To simplify the argument 
we  assume  $n = 1, \lambda  =  0$  and  $U_0(x)=  \delta^{(i)}(x)$.  
Representatives of $U$  resp. $\widetilde{U}$ are   
$u_\eps(x,t)=\mbox{Arsinh}(e^t \mbox{sinh}(\rho_\eps^{(i)}(x)))$ 
and $\widetilde{u}_\eps(x,t)=\mbox{Arsinh}(e^{c\eta+t}  
\mbox{sinh}(\rho_\eps^{(i)}(x)))$. For $\psi \in {\cal D}(\R^2)$ 
we have 
\begin{eqnarray*} 
& I_\eps^i := \int\!\int \widetilde{u}_\eps(x,t) \psi(x,t)dx dt= &\\ 
&=\int\!\int   \! \int_0^1  \theta(e^{c\eta+t},  \sigma 
\eps^{-i-1}\rho^{(i)}(x))   d\sigma   \eps^{-i}   \rho^{(i)}(x) 
\psi(\eps x,t) dx dt &  
\end{eqnarray*} 
where   $\theta(\alpha,y)   :=  \frac{d}{dy}\mbox{Arsinh}(\alpha 
\mbox{sinh}(y))$  for  $\alpha>0$, $y\in \R$. Since $\theta$ is 
bounded by  $\max(1,\alpha)$ and 
$\lim_{|y|\to\infty}\theta(\alpha,y)=1$    it    follows   that 
$I_\eps^0\to   \int   \psi(0,t)   dt$,  so  $\widetilde{U}$  is 
associated with the delta function on the $t$-axis, as desired. 
For $i\ge 1$ we write 
\begin{eqnarray*} 
& I_\eps^i =  
\int\!\int   \! \int_0^1  (\theta(e^{c\eta+t},  \sigma 
\eps^{-i-1}\rho^{(i)}(x)) -1)  d\sigma   \eps^{-i}   \rho^{(i)}(x) 
\psi(\eps   x,t)   dx  dt  + & \\ 
& + (-1)^i  \int\!\int  \rho(x)  \partial_x^i\psi(\eps x,t) dx dt &  
\end{eqnarray*} 
Here the  second  term converges to $(-1)^i\int \partial_x^i\psi(0,t)$ and 
the first term goes to zero since  
$\int_0^1     |\theta(\alpha,\sigma     y)    -1|d\sigma    \le 
\frac{2|\alpha^2-1|}{\alpha|y|}(1-e^{-|y|})$   for   $y\not=0$. 
This proves the claim for $\rho\in {\cal D}(\R)$. For  
$\rho\in {\cal S}(\R)$ splitting the $x$-integral into one from 
$-\frac{1}{\sqrt{\eps}}$  to  $\frac{1}{\sqrt{\eps}}$  and  one 
over  $|x|\ge  \frac{1}{\sqrt{\eps}}$  gives the same result. 
\et 
\subsection{Continuity Properties} \label{contprop} 
In this section we work out a different strategy for transferring classical 
point symmetries into the ${\cal G}$-setting. This approach, suggested 
in \cite{OR}, consists in a more topological way of looking at the 
transfer problem by using continuity properties of differential operators. 
As we have pointed out in the discussion following (\ref{differep}), the 
main obstacle against directly applying classical symmetry groups  
componentwise to representatives of generalized solutions is that the  
differential equations need not be satisfied componentwise. However, 
there are certain classes of partial differential operators that do  
allow such a direct application. Consider a linear partial differential 
operator $P$ giving rise to an equation 
\be \label{rightinv} 
PU = 0 
\end{equation} 
in ${\cal G}$ and let $G$ be a classical slowly increasing symmetry 
group of (\ref{rightinv}). Furthermore, suppose that $P$ possesses a 
continuous homogeneous (but not necessarily linear) right inverse  
$Q$. If $U=\cl[u]$ is 
a solution to (\ref{rightinv}) in ${\cal G}(\Omega)$ 
then there exists some $n\in {\cal N}(\Omega)$ such that 
\[ 
Pu = n. 
\] 
Since $Q$ is a right inverse of $P$ this implies 
\be \label{aa} 
P(u_\eps-Qn_\eps) = 0 \quad \forall \eps\in I. 
\end{equation} 
Also,  $Qn\in  {\cal  N}(\Omega)$  due  to  the  continuity and 
homogeneity assumption on $Q$. If $g\in G$, (\ref{aa}) implies 
\[ 
P(g(u_\eps - Qn_\eps)) = 0 \quad \forall \eps \in I. 
\]  
By definition, 
\[ 
P(gU) = \cl[P(gu)] = \cl[P(g(u-Qn))], 
\] 
so $gU$ is a solution as well. Summing up, $G$ is a symmetry group in  
${\cal G}$. The following result will serve to secure the existence of a  
right inverse as above for a large class of linear differential operators. 
\bp \label{michael} 
Let $E$, $F$ be Fr\'echet spaces and $A$ a continuous linear map from 
$E$ onto $F$. Then $A$ has a continuous homogeneous right inverse $B:F\to E$. 
\et

\pr See \cite{M}, p.~364.\ep\ms From these preparations we conclude 
\bt \label{cinv} 
Let  
\[ 
\Delta_{\nu}(x,u^{(n)})=0,\quad \nu=1,\ldots,l 
\] 
be a system of linear PDEs  
with slowly increasing $\Delta_{\nu}$ and let $\eta \rightarrow g_\eta$ 
be a slowly increasing symmetry group of this system. 
Assume that the operator defined by the left hand side is surjective 
$({\cal C}^\infty (\Omega))^l \rightarrow ({\cal C}^\infty (\Omega))^l\, .$ 
Then $\eta \rightarrow g_\eta$ is a symmetry group for the system in 
${\cal G}(\Omega)$ as well. \ep 
\et

The assumptions of Theorem \ref{cinv} are automatically satisfied for any 
linear partial differential operator with constant coefficients on an 
arbitrary convex open domain (see \cite{Hoe2}, 10.6). 
\bex  \label{acex} 
The system of one-dimensional linear acoustics 
\be \label{acoust} 
\begin{array}{c}    
P_t + U_x = 0\hphantom{.} \\ 
U_t + P_x = 0. 
\end{array} 
\end{equation} 
is transformed via $U = V - W, P = V + W$ into 
\be \label{transf}  
\begin{array}{c}    
\,V_t\, +\, V_x \,= 0\hphantom{.} \\ 
W_t - W_x = 0. 
\end{array} 
\end{equation} 
Using   the   infinitesimal   generators  $\Phi(v)\partial_v  + 
\Psi(w)\partial_w$  ($\Phi$, $\Psi$ arbitrary smooth functions) 
of   (\ref{transf})  we  obtain  symmetry  transformations  for 
(\ref{acoust}) of the form 
\begin{eqnarray*} 
\widetilde{U} = F^{-1} \left( \eta + F(\frac{1}{2}(P + U)) \right)  
    - G^{-1} \left( \theta + G(\frac{1}{2}(P - U)) \right) \\ 
\widetilde{P} = F^{-1} \left( \eta + F(\frac{1}{2}(P + U)) \right)  
    + G^{-1} \left( \theta + G(\frac{1}{2}(P - U)) \right)  
\end{eqnarray*} 
with arbitrary diffeomorphisms $F,G$. Since (\ref{acoust}) satisfies the  
assumptions  of  Theorem \ref{cinv}  on $\Omega = \R^2$ it follows that 
any slowly increasing transformation of this form is a symmetry 
of (\ref{acoust}). In particular, this includes nonlinear  
transformations  of distributional solutions, cf. Example \ref{nonlind}. 
\et 
In  the  remainder  of  this  section  we discuss the interplay 
between symmetry groups and solutions of PDEs in the sense of  
association. Consider 
\be \label{asseq} 
\Delta_\nu(x,u^{(n)}) \approx 0, \quad 1\le \nu \le l 
\end{equation} 
in ${\cal G}$. A slowly increasing symmetry group of the corresponding  
system 
\[  
\Delta(x,u^{(n)}) =  0, \quad 1\le \nu \le l 
\] 
is called a symmetry group in the sense of association if it 
transforms solutions of (\ref{asseq}) into other such solutions.  
The first question to be answered in this context is whether  
one can derive conditions on the  
form of the factorization (\ref{glg7}) that will yield symmetry 
groups  in  the  sense  of association. It is clear that a 
sufficient  condition  is  to suppose that $Q$ depends 
exclusively    on  $\eta$  and $x$. Distributional solutions to 
linear  PDEs  arise as a special case of (\ref{asseq}) and have 
been  treated  in  \cite{ga}.  There,  the validity of equation 
(\ref{glg7})  with  $Q$  depending  on  $\eta$  and $x$ only is 
actually used to {\em define} symmetry groups in ${\cal D}'$. 
In order to remain  within  the  classical  distributional   
framework,  the admissible group transformations in \cite{ga} are  
restricted to projectable ones acting linearly in   
the  dependent  variables.  On  the  other hand, the method 
developed   there   is  even  applicable  to  linear  equations 
containing  distributional terms which allows to use invariance 
methods to compute fundamental solutions. 
 
Second, if $u$ is a solution to $\Delta(x,u^{(n)})=0$ in  
${\cal G}(\Omega)$ possessing an associated distribution,  
one may ask for which group actions $g$ this implies  
that $gu$ as well possesses an associated distribution. This is 
certainly  the case for admissible transformations in the above 
sense. On the 
other  hand,  we  have already seen in Example \ref{grex0} that 
even genuinely nonlinear symmetry transformations may preserve  
association properties.  
 
The  next  example  shows  that  nonlinear  group  actions  may 
transform distributional solutions in Examples \ref{grex0} and  
\ref{acex}  into  more  complicated distributional solutions or 
into  generalized  solutions  in ${\cal G}(\R^2)$ not admitting 
associated distributions. 

\bex \label{nonlind} 
We  consider  the  equation  $U_t+\lambda  U_x  = 0$ arising in 
(\ref{semi}) with $n = 1$  or in (\ref{transf}). We have already 
observed that 
$\widetilde{U}   =   F^{-1}(\eta+F(U))$   defines   a  symmetry 
transformation  for arbitrary diffeomorphisms $F$. Here we take 
$F\in {\cal C}^\infty(\R)$, $F'>0$, $F(y)=\mbox{sign}(y) 
\sqrt{|y|}$  for $|y|\ge 1$. We wish to compute $\widetilde{U}$ 
when   $U\in   {\cal   G}(\R^2)$   is  a  delta  wave  solution 
$U(x,t)\approx  \delta^{(i)}(x-\lambda  t)$. We take $U$ as the 
class of  $\rho_\eps^{(i)}(x-\lambda t)$ with $\rho\in {\cal D} 
([-1,1])$. We have when $\eta\ge 0$: 
\begin{itemize} 
\item[(i)] If $i=0$, that is  $U\approx  \delta(x-\lambda  t)$, 
then $\widetilde{U}\approx F^{-1}(\eta+F(0)) + \delta(x-\lambda   
t)$; 
\item[(ii)] If $i=1$, that is  $U\approx \delta'(x-\lambda  t)$, 
then\\  
$\widetilde{U}\approx 
F^{-1}(\eta+F(0)) + 2\eta \int\sqrt{|\rho'(y)|}\, dy \, \, 
\delta(x-\lambda  t)+\delta'(x-\lambda t)$; 
\item[(iii)] If $i\ge 2$ then $\widetilde{U}$ does not admit an 
associated distribution. 
\end{itemize} 
To   see  this,  we  may  assume  that  $\lambda=0$  and  write 
$a_\eps(x):= \rho_\eps^{(i)}(x)$ for brevity.  
Note  that $F^{-1}(y) = \mbox{sign}(y)y^2$ for $|y|\ge 1$. 
Let    $A_\eps   =   \{x\in   [-\eps,\eps]   :   |a_\eps(x)|\ge 
(\eta+1)^2\}$. If $x\in A_\eps$ and $a_\eps(x)\ge 0$ then  
$\eta + F(a_\eps(x))\ge 1$ and $F^{-1}(\eta+F(a_\eps(x))) =  
\eta^2  + 2\eta \sqrt{a_\eps(x)} + a_\eps(x)$. Also, if 
$x\in A_\eps$ and $a_\eps(x)< 0$ then  
$\eta + F(a_\eps(x))\le -1$ and $F^{-1}(\eta+F(a_\eps(x))) =  
-\eta^2  + 2\eta \sqrt{|a_\eps(x)|} + a_\eps(x)$. The functions 
$F^{-1}(\eta+F(a_\eps))$, $|a_\eps(x)|$ and $\sqrt{|a_\eps(x)|}$ 
are bounded on the complement of $A_\eps$. Thus 
\begin{eqnarray*} 
&  \int\!\int_{-\eps}^\eps  F^{-1}(\eta + F(a_\eps(x))) \psi(x,t) 
dxdt = &\\ 
& =\int\! \int_{A_\eps}(\pm\eta^2  + 2\eta \sqrt{|a_\eps(x)|} + a_\eps(x)) 
\psi(x,t) dxdt +  O(\eps)   = &\\ 
& = \int \!\int_{-\eps}^\eps    (2\eta \sqrt{|a_\eps(x)|}   + 
a_\eps(x))\psi(x,t) dx dt + O(\eps)& 
\end{eqnarray*} 
while  
\[ 
\int\!\int_{|x|\ge \eps} F^{-1}(\eta +F(a_\eps(x)))\psi(x,t)\,dxdt \to 
F^{-1}(\eta + F(0)) \int\!\int \psi(x,t)\, dxdt 
\] 
It  follows  that  $F^{-1}(\eta  +  F(a_\eps(x)))$ converges in 
${\cal  D}'(\R^2)$ if and only if $2\eta \sqrt{|a_\eps|} + a_\eps$ 
admits  an associated distribution. A simple computation yields 
the particular results (i), (ii), (iii). 
\et 
\section{Generalized Group Actions}  \label{ggas} 
Although the methods introduced in the previous sections enable 
an application of large classes of classical symmetry groups to elements 
of Colombeau algebras, they are but the first step in a theory of  
generalized  group  analysis of differential equations. In this 
section  we  develop  an  extension  of  the  methods  of group 
analysis  that  will  allow  to  consider  symmetry  groups  of 
differential  equations whose actions are generalized functions 
themselves. 
\subsection{Generalized Transformation Groups} \label{gentransgr}
Simple examples indicate the necessity of extending the methods 
of  group  analysis  of PDEs to equations involving generalized 
functions themselves: 
\bex \label{gsgex} 
Considering  (\ref{semi})  in \gt with a {\em generalized} 
function $f = \cl[(f_\eps)_{\eps\in I}]$ $\in$ ${\cal G}_\tau$ we can 
apply  the  classical  algorithm  for calculating symmetry  
groups componentwise to the equations 
\[ 
\partial_t u_\eps + \lambda\cdot \nabla_x u_\eps = f_\eps(u_\eps) 
\] 
thereby  obtaining  infinitesimal  generators  with generalized 
coefficient  functions. Thus the question arises in which sense 
such generators induce symmetries of the differential equation. 
More generally, one can consider differential equations in  
${\cal G}_\tau$ of the form  
\[ 
P(x,U^{(n)})=0 
\] 
where $P$ is a generalized function. 
\et 
As is indicated by Example \ref{gsgex}, composition of generalized  
functions  will  inevitably  occur in a generalization of group 
analysis. For this purpose, we shall apply 
suitable  variants  of  Colombeau  algebras  for  the following 
considerations, namely ${\cal G}_\tau(\R^n)$ and  
$\widetilde{{\cal G}}_\tau(\R\times\R^n) = \widetilde{{\cal G}}_\tau 
(\R^{1+n})$.  
\bd \label{gga} 
A generalized group action  
on $\R^n$ is an element $\Phi$ of  
$(\widetilde{{\cal G}}_\tau(\R^{1+n}))^n$ such that: 
\begin{itemize} \rm 
\item[(i)] $\Phi(0,\,.\,)=\mbox{\rm id}$ in $({\cal G}_\tau(\R^{n}))^n$. 
\item[(ii)] $\Phi(\eta_1+\eta_2,\,.\,)=\Phi(\eta_1,\Phi(\eta_2,\,.\,))$ in 
$({\cal G}_\tau(\R^{2+n}))^n$. 
\end{itemize} 
\et
Before  we  turn to an infinitesimal description of generalized 
group actions let us shortly recall some basic definitions from 
\cite{point}  that are needed for a pointvalue characterization 
of generalized functions which in turn plays a fundamental role 
in the following considerations.  
Thus for any open set $\Omega\subseteq \R^n$ we set 
\[ 
\Omega_M := \{ (x_\eps)_{\eps\in I} \in \Omega^I:  
\exists p>0 \ \exists \ \eta>0 
\ |x_\eps|\le \eps^{-p} \ (0<\eps < \eta)\}. 
\] 
On $\Omega_M$ we define an equivalence relation by 
\[ 
(x_\eps)_{\eps\in I} \sim (y_\eps)_{\eps\in I} 
\ \Leftrightarrow \ \forall q>0 \ \exists \eta>0  
\ |x_\eps - y_\eps| \le \eps^q \ (0<\eps<\eta) 
\] 
and set $\widetilde{\Omega}:=\Omega_M/\sim$. $\widetilde{\Omega}$ 
is  called  the  set  of  generalized  points  corresponding to 
$\Omega$. The set of compactly supported points is defined as 
\[  
\widetilde{\Omega}_c = \{\widetilde{x}\in  
\widetilde{\Omega} : \exists 
\mbox{ representative } (x_\eps)_{\eps\in I} \  
\exists K\subset\subset \Omega 
\ \exists \eta>0 \,:\, x_\eps\in K,\, \eps\in (0,\eta)\}. 
\]  
Note  that  for $\Omega=\R$ we have $\widetilde{\Omega} = {\cal 
R}$. Theorems  2.4,  2.7  and  2.10  of  \cite{point} establish  
that elements of ${\cal G}(\Omega)$, $\widetilde{\cal G}_\tau(\Omega)$ 
or  $\widetilde{\cal G}_\tau(\Omega\times\Omega')$ are uniquely 
determined by their pointvalues in $\widetilde{\Omega}_c$, 
$\widetilde{\Omega}$, or $\widetilde{\Omega}_c\times  
\widetilde{\Omega}'$, respectively.  For the theory of ODEs in 
the Colombeau framework we refer to \cite{HO}. 
\bd Let $\xi=(\xi_1,\ldots,\xi_n)\in ({\cal G}_\tau(\R^n))^n$. 
The generalized vector field  $X=\sum\limits_{i=1}^n \xi_i(x)  
\partial_{x_i}$ is called ${\cal G}$-complete if the initial  
value problem 
\beast 
&& \dot{x}(t) = \xi(x(t)) \\ 
&& x(t_o) = \widetilde{x}_o 
\eeast 
is uniquely solvable in ${\cal G}(\R)^n$ for any  
$\widetilde{x}_o\in {\cal R}^n$ and any $t_o\in \R$. 
\et 
\bd 
Let $\Phi$ be a generalized group action on $\R^n$ and set  
\[ 
\xi:=\deta\Phi(\eta,\,.\,) \, \in ({\cal G}_\tau(\R^n))^n.  
\] 
If the generalized vector field 
$X=\sum\limits_{i=1}^n \xi_i(x)\partial_{x_i}$ is ${\cal G}$-complete,  
then $X$ is called the infinitesimal generator of $\Phi$.  
In this case, $\Phi$ is also called ${\cal G}$-complete. 
\et
By \cite{HO}, every generalized vector field with ${\cal G}_\tau$-components  
whose gradient is of $L^\infty$-log-type is ${\cal G}$-complete.  
The notion of infinitesimal generator is well-defined due to 
\bp \label{gen}  
Every ${\cal G}$-complete generalized group action is uniquely  
determined by its infinitesimal generator. 
\et  

\pr Let $\Phi'$, $\Phi''$ be two \gc generalized group actions  
with the same infinitesimal  
generator $X=\sum_{i=1}^n \xi_i(x)\partial_{x_i}$. Then  
both functions satisfy 
\[ 
\frac{d}{d\eta}\Phi(\eta,x)=\frac{d}{d\mu}{\Big\vert}_{_{0}} 
\Phi(\eta+\mu,x)= 
\frac{d}{d\mu}{\Big\vert}_{_{0}}\Phi(\mu,\Phi(\eta,x))=\xi(\Phi(\eta,x)). 
\] 
Now given any $\widetilde{x}\in {\cal R}^n$, it follows that both 
$\eta \to\Phi'(\eta,\widetilde{x})$ and $\eta \to \Phi''(\eta,\widetilde{x})$ 
solve the initial value problem 
\beast 
&& \dot{x}(\eta) = \xi(x(\eta)) \\ 
&& x(0) = \widetilde{x} 
\eeast 
By assumption this entails that $\Phi'(\,.\,,\widetilde{x})= 
\Phi''(\,.\,,\widetilde{x})$ in $({\cal G}(\R))^n$. 
Consequently, 
\[ 
\Phi'(\widetilde{\eta},\widetilde{x}) =  
\Phi''(\widetilde{\eta},\widetilde{x}) 
\] 
for all $\widetilde{\eta}\in {\cal R}_c$ and all  
$\widetilde{x}\in {\cal R}^n$. 
The claim now follows from \cite{point}, Theorem 2.10.\ep\ms 
As in the classical theory, we are first going to  
investigate symmetry groups of algebraic equations: 
\bd 
Let $F\in {\cal G}_\tau (\R^n)$ and let $\Phi$ be a generalized 
group action on $\R^n$. $\Phi$ is called a symmetry group of  
the equation  
\[ 
F(x)=0 
\] 
in \gtrn if for any $\widetilde{x}\in {\cal R}^n$  
with $F(\widetilde{x})=0\in {\cal R}$ 
it follows that $\eta\rightarrow F(\Phi(\eta,\widetilde{x}))=0$  
in ${\cal G}(\R)$ (or, equivalently, $F(\Phi(\widetilde{\eta}, 
\widetilde{x}))=0$ in 
${\cal R}$ for every $\widetilde{\eta} \in {\cal R}_c$). 
\et 
A characterization of symmetry groups of (generalized) algebraic equations 
in terms of infinitesimal generators is provided by 
\bt \label{algsymg}  
Let $F\in\gtn$ be of the form 
\[ 
F(x_1,\ldots,x_n)=x_i-f(x_1,\ldots,x_{i-1},x_{i+1},\ldots,x_n) 
\] 
for some $1\le i\le n$ and $f\in {\cal G}_\tau(\R^{n-1})$. Let $\Phi$ be 
a ${\cal G}$-complete generalized group action with infinitesimal generator  
$X=\sum_{i=1}^n \xi_i(x)\partial_{x_i}$ and suppose that $x'\rightarrow 
\xi(x',f(x'))$ defines a generalized vector field on $\R^{n-1}$ such that  
the corresponding system of ODEs possesses a flow in 
$(\widetilde{{\cal G}}_\tau(\R^{1+(n-1)}))^{n-1}$.  
The following conditions are equivalent: 
\begin{itemize} 
\item[(i)] $\Phi$ is a symmetry group of $F(x)=0$. 
\item[(ii)] If $\widetilde{x}\in  {\cal R}^n$ with  
$F(\widetilde{x})=0\in {\cal R}$ 
it follows that $X(F)(\widetilde{x})=0$ in ${\cal R}$. 
\end{itemize} 
\et 

\pr (i) $\Rightarrow$ (ii): Consider the function $(\eta,x)\rightarrow 
F(\Phi(\eta,x))\in \widetilde{\cal G}_\tau(\R^{1+n})$. We have 
\[\frac{d}{d\eta} F(\Phi(\eta,x))=\sum_{i=1}^n 
\frac{\partial F}{\partial x_i} 
(\Phi(\eta,x))\xi_i(\Phi(\eta,x))=X(F)(\Phi(\eta,x)),\] 
so that $\detas F(\Phi(\eta,x))=X(F)(x)$ in $\gtn$. 
Let $\widetilde{x}\in{\cal R}^n$ such that $F(\widetilde{x})=0$. Then 
$F(\Phi(\,.\,,\widetilde{x}))=0$ in ${\cal G}(\R)$. 
Thus $\detas F(\Phi(\eta,\widetilde{x}))=0$ in  
${\cal R}$ which means that $X(F)(\widetilde{x})=0$ in 
${\cal R}$.\\ 
(ii) $\Rightarrow$ (i):  
We assume $F(x_1,\ldots,x_n)=x_n-f(x_1,\ldots,x_{n-1})$ 
and abbreviate $(x_1,$ $\ldots$, $x_{n-1})$ by $x'$.  
Our first claim is that 
\[ 
\xi_n(x',f(x'))=\sum_{j=1}^{n-1}\xi_j(x',f(x'))\partial_j f(x') 
\mbox{  in  } {\cal G}_\tau(\R^{n-1})  
\] 
Indeed, if $\widetilde{x}'\in {\cal R}^{n-1}$ then  
$F(\widetilde{x}',f(\widetilde{x}'))=0$ in ${\cal R}$.  
Hence $X(F)(\widetilde{x}',f(\widetilde{x}'))=0$ in ${\cal R}$ for 
all $\widetilde{x}'$ by our assumption. Our claim now follows from  
\cite{point}, Theorem 2.7. Consider the following system of ODEs in  
${\cal G}_\tau$: 
\beast 
   \dot{x}_j(t) &=& \xi_j(x',f(x')) 
   \quad (j=1,\ldots,n-1)\\ 
   x'(0) &=& \widetilde{a}'\in {\cal R}^{n-1} 
\eeast 
By our assumption, this system has a flow $(\eta,a')$ $\rightarrow$ 
$(h_1(\eta,a'),$ $\ldots,$ $h_{n-1}(\eta,a'))$ 
in $(\widetilde{{\cal G}}_\tau(\R^{1+(n-1)}))^{n-1}$. Set $g_n(\eta,a):=$ 
$f(h_1(\eta,a')$,\ldots,$h_{n-1}(\eta,a'))$. Then $g_n(0,a)=f(a')$ and  
\[ 
g(\eta,a)=(g_1(\eta,a),\ldots,g_n(\eta,a)):=  
(h_1(\eta,a'),\ldots,h_{n-1}(\eta,a'),g_n(\eta,a)) 
\] 
is in $(\widetilde{{\cal G}}_\tau(\R^{1+n}))^{n}$. 
If $\widetilde{a}\in {\cal R}^n$ then $F(g(\eta,\widetilde{a}))=0$  
in ${\cal R}$ for all $\eta \in {\cal R}_c$. Therefore, if we can show that  
$g(\,.\,,\widetilde{a})=\Phi(\,.\,,\widetilde{a})$ in $({\cal G}(\R))^n$  
for all $\widetilde{a}$ with $F(\widetilde{a})=0$,  
the proof is completed. Now we have $\dot{g}_j(\eta,a)=\xi_j( 
g_1(\eta,a),\ldots,g_n(\eta,a))$ for $1\le j\le n-1$ and 
\beast 
& \dot{g}_n(\eta,a)=\sum\limits_{i=1}^{n-1}\frac{\partial f}{\partial x_i} 
(g_1(\eta,a),\ldots,g_{n-1}(\eta,a))\dot{g_i}(\eta,a)= &\\ 
& =\xi_n(g_1(\eta,a),\ldots,f(g_1(\eta,a),\ldots,g_{n-1}(\eta,a)))= 
\xi_n(g(\eta,a)).& 
\eeast 
If $F(\widetilde{a})=0$ in ${\cal R}$ then $\widetilde{a}_n= 
f(\widetilde{a}')$, so that $g(0,\widetilde{a})=(\widetilde{a}', 
f(\widetilde{a}'))=\widetilde{a}=\Phi(0,\widetilde{a})$. 
Thus $g(\,.\,,\widetilde{a})$ and $\Phi(\,.\,,\widetilde{a})$ solve the same  
initial value problem. Since $X$ is ${\cal G}$-complete,  
the claim follows.\ep
\subsection{Symmetries of Differential Equations} \label{symdiffeq}
In this section we are going to apply the above results to  
symmetry groups of differential equations involving generalized 
functions. To this end, we will 
first have to define generalized group actions on generalized functions.  
Once we have done this, by a symmetry group of a differential equation 
we will again mean a group action that transforms solutions into other  
solutions. 
Thus, from now on we will exclusively consider group actions on some space  
$\R^p\times \R^q$ of independent and dependent variables. 
\bd 
A generalized group action $\Phi\in  
(\widetilde{{\cal G}}_\tau (\R\times \R^{p+q}))^{p+q}$ is called  
projectable if it is of the form  
\[ 
\Phi(\eta,(x,u))=(\Xi_\eta(x),\Psi_\eta(x,u)), 
\] 
where $\Xi\in (\widetilde{{\cal G}}_\tau (\R\times\R^p))^p$ and 
$\Psi\in (\widetilde{{\cal G}}_\tau (\R\times \R^{p+q}))^q$. 
\et 
The group properties in this case read: 
\bea 
& \Xi_{\eta_1+\eta_2} = \Xi_{\eta_1}\circ \Xi_{\eta_2}  
\quad \mbox{in } {\cal G}_\tau(\R^p)\ \forall \eta_1,\eta_2 \in {\cal R}_c.
\label{gr1} & \\ 
& \Psi_{\eta_1+\eta_2} (x,u) = \Psi_{\eta_1}(\Xi_{\eta_2}(x),\Psi_{\eta_2} 
(x,u)) \quad \mbox{in } {\cal G}_\tau(\R^{p+q})\ \forall \eta_1,\eta_2  
\in {\cal R}_c. \label{gr2} & 
\eea 
In particular, we have  
\be\label{equa} 
\Xi_\eta\circ\Xi_{-\eta}=\,\mbox{id}\ \quad \mbox{in }  
{\cal G}_\tau (\R^p)\ \forall\eta\in {\cal R}_c.  
\end{equation} 
An adaptation of Lie group analysis to spaces of distributions faces  
the fundamental problem that while the methods of classical Lie 
group analysis of differential equations are {\em geometric} in the  
sense that group action on functions is defined via graphs,  
in classical distribution theory there is no means of defining graphs 
of    distributions.    However,    due   to   the   pointvalue 
characterization obtained in \cite{point} this problem can be dealt  
with in a satisfactory manner within Colombeau algebras:  
\bd
Let $U\in ({\cal G}(\R^p))^q$ and $V\in ({\cal G}_\tau(\R^p))^q$. The graphs 
of $U$ and $V$ are defined as 
\beast 
& \Gamma_U:=\{(\widetilde{x},U(\widetilde{x}))\mbox{ : } 
\widetilde{x}\in {\cal R}^p_c\}\hphantom{.} &\\ 
& \Gamma_V:=\{(\widetilde{x},U(\widetilde{x}))\mbox{ : } 
\widetilde{x}\in {\cal R}^p\}. &  
\eeast 
\et
It follows directly from \cite{point}, Theorems 2.4  and 2.7 that any  
generalized function is uniquely determined by its graph.  
Our next aim is to define generalized group actions  
on generalized functions. As in the classical case this is done  
geometrically,  i.e.~by transformation of graphs. The following 
result is immediate from the definitions: 
\bp  
Let $U\in ({\cal G}_\tau(\R^p))^q$ and let $\Phi$ be a  
projectable generalized group action on $\R^p\times \R^q$. Then  
$\Phi_\eta(\Gamma_U)=\Gamma_{\Phi_\eta(U)}$  
in ${\cal R}^{p+q}$ for each $\eta$, where $\Phi_\eta(U)$ denotes the  
element 
\[ 
x\rightarrow \Psi_\eta(\Xi_{-\eta}(x),U\circ \Xi_{-\eta}(x)) 
\] 
of $({\cal G}_\tau(\R^p))^q$.\ep 
\et

We are now able to give a geometric characterization of solutions of PDEs 
in ${\cal G}_\tau$. 
\bp \label{algeq}  
Consider the system of PDEs  
\be \label{dsystg} 
\Delta_\nu (x,U^{(n)}) = 0  \quad 1\le \nu \le l 
\end{equation} 
in ${\cal G}_\tau(\R^p))^q$ (where 
$\Delta\in ({\cal G}_\tau((\R^p\times \R^q)^{(n)}))^l$). Set  
\[ 
{\cal S}_\Delta := 
\{\widetilde{z}\in {\cal R}^{(n)}\, : \, \Delta_\nu(\widetilde{z})=0  
\ (1\le \nu \le l)\}. 
\] 
Then $U\in ({\cal G}_\tau(\R^p))^q$ is a 
solution of the system iff $\Gamma_{\prol U}\subseteq {\cal S}_\Delta$. 
\et 

\pr This follows immediately from \cite{point}, Theorem 2.7. \ep\ms 
Prolongation of generalized group actions can be handled 
in a similar fashion as in the classical theory. 
Thus, let $\Phi$ be a projectable generalized group action on  
$\R^p\times\R^q$. We want to define  
the $n$-th prolongation $\prol \Phi$ as a projectable generalized group  
action on $(\R^p\times \R^q)^{(n)}$. Let $z\in  
(\R^p\times \R^q)^{(n)}$  
and choose $h\in {\cal O}_M(\R^p)^q$ such  
that $(z_1,\ldots,z_p,\prol h(z_1,\ldots,z_p))=z$. Now set  
\be \label{proldefg} 
\prol \Phi(\eta,z) := 
(\Xi_{\eta}(z_1,...,z_p),\prol(\Phi_\eta(h))(\Xi_{\eta}(z_1,...,z_p))). 
\end{equation} 
Using for $h$ a suitable Taylor polynomial, it follows that 
$\prol \Phi\in (\widetilde{{\cal G}}_\tau (\R\times(\R^p\times \R^q)^ 
{(n)})^N$ (where $N=\mbox{dim}((\R^{p+q})^{(n)})$). Moreover, the definition  
does not depend on the particular choice of $h$, which follows exactly 
as in the classical case.  
\blem \label{save} 
Let $\widetilde{z}\in (\gR^p \times \gR^q)^{(n)}$ and assume that $U\in  
({\cal G}_\tau(\R^p))^q$ satisfies  
$(\widetilde{z}_1,$ $\ldots,$ $\widetilde{z}_p,$ $\prol U  
(\widetilde{z}_1,\ldots, 
\widetilde{z}_p)) = \widetilde{z}$. Then 
\be \label{proldefpv} 
\prol \Phi (\eta,\widetilde{z}) = (\Xi_\eta(\widetilde{z}_1,\ldots, 
\widetilde{z}_p), \prol (\Phi_\eta(U))(\Xi_\eta(\widetilde{z}_1,\ldots, 
\widetilde{z}_p)))\quad \forall \eta \in {\cal R}_c. 
\end{equation} 
\et 

\pr Let $U=\cl[(u_\eps)_{\eps\in I}]$ and choose a 
representative $(z_\eps)_{\eps\in I}$ 
of $\widetilde{z}$ 
such that  
\[ 
(z_{1\eps},\ldots,z_{p\eps},\prol u_\eps(z_{1\eps},\ldots,z_{p\eps})) = z_\eps 
\quad \forall \eps. 
\] 
Using the chain rule as in Proposition 
\ref{prolprop}, it follows that the right hand sides of (\ref{proldefg}) 
(with $z$ replaced by $\widetilde{z}$) and of (\ref{proldefpv}) have  
the same representative (depending {\em exclusively} on 
$(z_\eps)_{\eps\in I}$).\ep 
\bp 
$\prol \Phi$ is a generalized group action on $(\R^p\times \R^q)^{(n)}$. 
\et

\pr  Property \ref{gga}  (i)  is  clearly satisfied.  Concerning  (ii),  
according to \cite{point}, Theorem 2.7 it suffices to show that 
\[ 
\prol \Phi(\eta_1+\eta_2, \widetilde{z}) = \prol \Phi(\eta_1,\prol\Phi 
(\eta_2,\widetilde{z}))\quad \forall \eta_1,\eta_2\in {\cal R}_c, \ \forall 
\widetilde{z} \in (\gR^p\times \gR^q)^{(n)}. 
\] 
Choose some $U\in ({\cal G}_\tau(\R^p))^q$ with  
$(\widetilde{z}_1,\ldots,\widetilde{z}_p,\prol U  
(\widetilde{z}_1,\ldots,\widetilde{z}_p))   =   \widetilde{z}$. 
Then due to Lemma \ref{save} we have 
\[ 
\prol \Phi (\eta_2,\widetilde{z}) = (\Xi_{\eta_2}(\widetilde{z}_1,\ldots, 
\widetilde{z}_p), \prol (\Phi_{\eta_2}(U))(\Xi_{\eta_2}( 
\widetilde{z}_1,\ldots,\widetilde{z}_p))). 
\] 
By (\ref{proldefpv}) this implies 
$\prol \Phi(\eta_1, \prol \Phi(\eta_2, \widetilde{z})) =  
\prol \Phi(\eta_1+\eta_2, \widetilde{z})$.  
\ep\ms 
As in the classical case we therefore have (using the notations from 
Proposition \ref{algeq}): 
\bp\label{prolg}  
Let $\Phi$ be a projectable generalized group action  
on $\R^p\times \R^q$ such that $\prol\Phi$ is a  
symmetry group of the algebraic equation 
$\Delta(z)=0$. Then $\Phi$ is a symmetry group of (\ref{dsystg}). 
\et

\pr If $U\in {\cal G}_\tau(\R^p)$ is a solution of (\ref{dsystg}) then     
$\Gamma_{\prol U}\subseteq {\cal S}_\Delta$ by Proposition \ref{algeq}. Thus  
\[ 
\Gamma_{\prol (\Phi_\eta U)}= 
\prol\Phi_\eta (\Gamma_{\prol U})\subseteq {\cal S}_\Delta,  
\] 
so that, again from Proposition \ref{algeq}, the claim follows. \ep 
\bd Let $X$ be a ${\cal G}$-complete generalized vector field.  
The $n$-th prolongation  
of $X$ is defined as the infinitesimal generator of the $n$-th  
prolongation of the generalized group action $\Phi$ corresponding to $X$: 
\[ 
\prol X|_z = \deta \prol \Phi_\eta(z), 
\] 
provided that $\prol \Phi$ is ${\cal G}$-complete as well.  
In this case, both $X$ and $\Phi$ are called  
${\cal G}$-$n$-complete. 
\et \par From Theorem \ref{algsymg} and Proposition \ref{prolg} 
we immediately conclude 
\bt \label{sgth}  
Under the assumptions of Proposition \ref{algeq}, let $\Phi$ be a  
${\cal G}$-$n$-complete generalized group action on $\R^p\times \R^q$ 
with  infinitesimal  generator  $X$ such that the conditions of 
Theorem \ref{algsymg} are satisfied for $\Delta$ and $\prol \Phi$.  
If  
\[ 
\prol X(\Delta)(\widetilde{z})=0 
\quad \forall \widetilde{z}\in ({\cal R}^p\times {\cal R}^q)^{(n)} 
\mbox{ with } \Delta(\widetilde{z})=0, 
\]  
then $\Phi$ is a symmetry group of (\ref{dsystg}).\ep 
\et 

In order to be able to apply the same algorithm as in classical Lie theory 
for the determination of the symmetry group of a generalized 
PDE, the final step is to verify that the formulas for prolongation of vector 
fields carry over to generalized vector fields. 
\bt\label{prolform}  
Let  
\[ 
X=(x,u)\rightarrow\sum_{i=1}^p \xi_i(x)\partial_{x_i} + \sum_{\alpha=1}^q 
\psi_\alpha(x,u)\partial_{u^\alpha} 
\] 
be a ${\cal G}$-$n$-complete generalized vector field with  
corresponding projectable group action $\Phi$ on $(\R^p\times \R^q)$. Then 
\[ 
\prol X=X+\sum_{\alpha=1}^q\sum_J \psi_\alpha^J(x,u^{(n)})\partial_{u_J 
^\alpha} 
\] 
where $J=(j_1,...,j_k)$, $1\le j_k\le p$ for $1\le k\le n$ and 
\[ 
\psi_\alpha^J(x,u^{(n)})=D_J(\psi_\alpha - \sum_{i=1}^p \xi_i u_i^\alpha) 
+ \sum_{i=1}^p \xi_i u_{J,i}^\alpha 
\] 
\et 

\pr  Using  the  machinery  developed  so  far, this is an easy 
modification   of  the  proof  of  the  classical  result  (see 
\cite{Olv}, Theorem 2.36). \ep\ms 
We may summarize the results of this section as follows:   
In  order  to determine the symmetries of a differential 
equation  involving generalized functions, the algorithm (as in 
the  classical case) is to make an ansatz for the infinitesimal 
generators,    calculate   their   prolongations   according to 
Theorem \ref{prolform}   and  then  use  Theorem \ref{sgth}  to 
determine the defining equations for the coefficient  functions  
of the infinitesimal generators. The defining equations now  
yield PDEs in ${\cal G}_\tau$. Any solution 
of these equations that defines a ${\cal G}$-$n$-complete 
generator will upon integration yield a symmetry group in  
${\cal G}_\tau$. 
\bex \label{gengrex} 
Scalar conservation laws of the form
\be \label{gtcons} 
u_t + F(u) u_x = 0 
\end{equation} 
arise in the kinetic theory of traffic flow. Here $u$ denotes the
density, and the propagation velocity $F$ may be a strictly
decreasing function of $u$ with one or more jumps. A typical case
is a unimodal flux function (whose derivative is $F$) with a kink at
its maximum, as supported by experimental data \cite{Hall}.
Convolution with a nonnegative mollifier $(\rho_\eps)_{\eps\in I}$
allows to interpret $F$ as an element of  ${\cal  G}_\tau(\R)$
which is invertible. Thus our theory of symmetry transformations for
equations with generalized nonlinearities applies. The 
determining equations are
\begin{eqnarray*} 
& \varphi_t + F \varphi_x = 0  & \\ 
& -\xi_x  + F \tau_t + \tau F_t + \varphi F_u - F \xi_x + F^2 \tau_x 
+ \xi F_x = 0 & 
\end{eqnarray*} 
with  infinitesimal  generator  ${\bf v} = \xi(x,t)\partial_x + 
\tau(x,t)\partial_t + \varphi(x,t,u)\partial_u$. 
As a particular solution we obtain ${\bf v} = xt\partial_x + 
t^2\partial_t    +    (F'(u))^{-1} (x - t F(u))\partial_u$. The 
corresponding   generalized  group  action  can  be  calculated 
explicitly  in  \gt  showing  that  if $u$ is a \gt-solution to 
(\ref{gtcons}) then so is  
\[ 
(x,t)  \to  F^{-1}\left(\eta  x  (1+\eta  t)^{-1}  +  F(u(x  (1+\eta 
t)^{-1},t (1+\eta t)^{-1}) (1+\eta t)^{-1}\right) 
\] 
In particular, a constant state $u$ is transformed into a generalized
solution to (\ref{gtcons}) which, depending on the shape of $F$,
will generally be associated with a piecewise smooth function.
\et
\bex \label{dAlHam} 
The nonlinear d'Alembert-Hamilton system
\be \label{dAlsys} 
\begin{array}{l} 
u_{tt} - u_{xx} - u_{yy} - u_{zz}  = F(u) \\   
u_t^2 - u_x^2 - u_y^2 - u_z^2 = G(u)
\end{array} 
\end{equation} 
arises in the study of relativistic field equations \cite{Cie} and as
a constraint in reducing the nonlinear wave equation to an ODE 
\cite{Fush1, Fush2}. One of its symmetries is generated by the vector field
${\bf v} = \varphi (u)\partial_u$ where the function $\varphi$ has to
satisfy
\begin{eqnarray*} 
& F\varphi_u - \varphi F_u + G\varphi_{uu} = 0  & \\ 
& 2G\varphi_u - \varphi G_u = 0\,. & 
\end{eqnarray*} 
In particular, in the isotropic case $F \equiv G \equiv 0$ the function
$\varphi$ is arbitrary. In our theory it may be taken in 
${\cal  G}_\tau(\R)$ subject to the ${\cal G}$-completeness conditions
formulated above.
As an example of the possible behavior of generalized transformations,
consider the vectorfield ${\bf v} = \varphi (u)\partial_u$ where
$\varphi \in {\cal  G}_\tau(\R)$ is the class of $(\varphi_\eps)_{\eps\in I}$
with $\varphi_\eps (u) = \tanh (\frac{u}{\eps})$. Thus $\varphi (u)$ is
associated with the jump function $- \mbox{sgn}  (u)$. Starting with a
classical smooth solution $u = u(x,t) \in {\cal O}_C(\R^4)$ of the
isotropic d'Alembert-Hamilton system ((\ref{dAlsys}) with $F \equiv
G \equiv 0$), the generalized symmetry transform generated by the
vector field ${\bf v}$ turns $u(x,t)$ into the generalized solution
$\tilde{U} \in {\cal G}_\tau(\R^4)$ with representative
\[
   \tilde{u}_\eps (x,t) = 
      \eps\,\mbox{Arsinh}\left( e^{\eta/\eps}\sinh\frac{u(x,t)}{\eps}\right)\,.
\]
When $\eta > 0$, it is straightforward to check that $\tilde{U}$ is associated
with the piecewise smooth function $v(x,t) = u(x,t) + \eta\, \mbox{sgn}  
(u(x,t))$. The generalized symmetry this way transforms smooth solutions into
discontinuous solutions.
\et  

{\bf  Acknowledgements.}  We  would  like  to thank M.~Grosser, 
G.~H\"ormann and P.~J.~Olver for several helpful discussions. A number
of constructive suggestions of the two referees led to improvements 
in the paper.

{\small
Michael Kunzinger: Universit\"at Wien, Institut f\"ur Mathematik, Strudlhofg. 4, 
A-1090 Wien, AUSTRIA;\\ {\it Email:} Michael.Kunzinger@univie.ac.at

Michael Oberguggenberger: Universit\"at Innsbruck, Institut f\"ur Mathematik und Geometrie,
Technikerstr. 13, A-6020 Innsbruck, AUSTRIA;\\ {\it Email:} michael@mat1.uibk.ac.at
}

\end{document}